\documentclass[11pt]{article}
\usepackage{amsmath,amsthm,amssymb}

\newtheorem{theorem}{Theorem}
\newtheorem{lemma}{Lemma}

\newtheorem{conjecture}{Conjecture}

\newcommand{\hee}{{\mathcal H}}
\newcommand{\fee}{{\mathcal F}}

\title{A Note on $G$-intersecting Families} \author{Tom
Bohman\thanks{Supported in part by NSF grant DMS-0100400.}
\hspace{.5in} Ryan R. Martin\thanks{Supported in part by NSF VIGRE
Grant DMS-9819950} \\ \\ Department of Mathematical Sciences \\ Carnegie
Mellon University \\ Pittsburgh, PA 15213}\date{}

\begin{document}

\maketitle 

\begin{abstract} 
Consider a graph $G$ and a \(k\)-uniform 
hypergraph $\hee$ on
common vertex set $[n]$.  We say that $\hee$ is 
{\it $G$-intersecting} if for every 
pair of edges in \(X,Y \in \hee \) there 
are vertices \(x \in X \) and \( y \in Y \) such that
\(x = y \) or \(x \) and \(y\) are joined by an 
edge in \(G\).  
This notion was introduced by Bohman, Frieze, 
Ruszink\'o and
Thoma who proved a natural generalization 
of the Erd\H{o}s-Ko-Rado Theorem for \(G\)-intersecting
\(k\)-uniform hypergraphs for \(G\) sparse and
\( k = O( n^{1/4} ) \).
In this note, we extend this result to 
\( k = O\left( \sqrt{n} \right) \).
\end{abstract}

\section{Introduction}

A hypergraph is said to be {\it intersecting} if every pair of edges has a 
nonempty intersection.  The well-known theorem of Erd\H{o}s, Ko and 
Rado~\cite{DF,EKR} details the extremal $k$-uniform intersecting hypergraph 
on $n$ vertices. 
\begin{theorem}[Erd\H{o}s-Ko-Rado] Let $k \le n/2$ and 
$\cal H$ be a $k$-uniform, intersecting hypergraph on vertex set $[n]$.   
We have \( |\hee| \leq \binom{n-1}{k-1} \).
Furthermore, $ |\hee| = \binom{n-1}{k-1} $ if 
and only if there exists \(v \in [n]\) such that 
\( \hee = \{ e \in \binom{n}{k} : v \in e \}\).
\label{tEKR}
\end{theorem}
\noindent
Of course, for \(k >n/2\) the hypergraph consisting of all \(k\)-sets
is intersecting.  So, extremal \(k\)-intersecting hypergraphs come in
one of two forms, depending on the value of \(k\). \\

Bohman, Frieze, Ruszink\'o and Thoma~\cite{BFRT} 
introduced a generalization of the notion of an 
intersecting hypergraph.  Let $G$ be a 
graph on a vertex set $[n]$ and $\cal H$ be a hypergraph, also on vertex set 
$[n]$.  We say {\it $\cal H$ is $G$-intersecting} if for any $e,f\in {\cal 
H}$, we have $e\cap f\neq\emptyset$ or there are vertices $v,w$ with $v\in 
e$, $w\in f$ and $v\sim_G w$.  We are intersected in the size and
structure of maximum \(G\)-intersecting hypergraphs; in particular,
we investigate
\[ N(G,k)=\max\left\{|{\cal H}| : {\cal H}\subseteq {[n]\choose k}
   \mbox{ and $\cal H$ is $G$-intersecting}\right\} . \] 
Clearly, Erd\H{o}s-Ko-Rado gives the
value of $N(E_n,k)$ where \(E_n\) is the empty graph 
on vertex set
\([n]\).  For a discusssion of \( N(G,k) \) 
for some other specific graphs see \cite{BFRT}. \\

In this note we restrict our attention to sparse graphs: 
those graphs for which \(n\) is large and the maximum degree 
of $G$, $\Delta(G)$, is a constant in \(n\).  What form can a 
maximum \(G\)-intersecting family take?  
If \(K\) is a maximum clique in \(G\) then a candidate 
for a maximum \(G\)-intersecting family is
\begin{equation*}
\hee_K := \left\{ X \in \binom{[n]}{k} : X \cap K 
\neq \emptyset \right\}.
\end{equation*} 
Note that such a hypergraph 
can be viewed as a natural generalization of the maximum 
intersecting hypergraphs given by 
Erd\H{o}s-Ko-Rado.  However, for many graphs and 
maximum cliques \(K\) 
one can add hyperedges to \( \hee_K\) to obtain a 
larger \(G\)-intersecting hypergraph. \\

Consider, for example,
\(C_n\), the cycle on vertex set \([n]\) (i.e. the graph on
\([n]\) in which \(u\) and \(v\) are adjacent iff 
\( u-v \in \{1,n-1\} \; {\rm  mod } \; n\)).  The set 
\( \{2,3\} \) is a maximum
clique in \(C_n\) and the set
\begin{equation}
\label{eq:set}
 \hee_{ \{2,3\}} \cup \left\{ X \in \binom{[n]}{k} : 
\{1,4\} \subseteq X  \right\} 
\end{equation}
is \(G\) intersecting.  Bohman, Frieze, Ruszink\'o and 
Thoma showed that 
\begin{equation}
\label{eq:cycle}
N(C_n,k) = \binom{n}{k} - \binom{n-2}{k} 
+ \binom{n-4}{k-2} 
\end{equation}
(i.e. the hypergraph given in (\ref{eq:set}) is maximum)
for \( k \) less than a certain constant times \( n^{1/4} \).
In fact, they showed that for arbitrary sparse graphs 
and \(k\) small, 
\( N(G,k) \) is given by a hypergraph that consists of 
\( \hee_K \) for some clique \(K\) together with a number of
`extra' hyperedges that cover the clique \(K\) in \(G\) (see
Theorem~1 of \cite{BFRT}).  
In this note we extend this result to larger values of \(k\). \\

\begin{theorem}
\label{tSQRT}
Let \(G\) be a graph on \(n\) vertices with maximum 
degree $\Delta$ and clique number $\omega$.  
There exists a constant $C$ (depending only on 
\(\Delta\) and \(\omega\)) such that if \( {\cal H} \) is a
\(G\)-intersecting \(k\)-uniform hypergraph and
\( k < C n^{1/2} \) then
   \[ |{\cal H}|\leq {n\choose k}-{n-\omega\choose k}
   + {\omega(\Delta-\omega + 1)\choose
   2}{n-\omega-2\choose k-2}. \]
   Furthermore, if  \( {\cal H} \) is a \(G\)-intersecting family of
   maximum cardinality
   then there exists a maximum clique \(K\) in \(G\) such that \( {\cal H} \) 
   contains
   all $k$-sets that intersect \(K\).
\end{theorem}
\noindent
An immediate corollary of this Theorem is that (\ref{eq:cycle}) holds
for \( k < C \sqrt{n} \). \\

Of course, a maximum \(G\)-intersecting hypergraph will not be of
the form `\( \hee_K \) together with some extra hyperedges' if \(k\)
is too large.  Even for sparse graphs, when \(k\) is large enough, 
there are hypergraphs that consist of nearly all of 
\( \binom{[n]}{k} \) that are \(G\)-intersecting.  In particular, 
Bohman, Frieze, 
Ruszink\'o and Thoma showed that if 
\(G\) is a sparse graph with minimum degree 
\( \delta \), \(c\) is a constant such that 
\( c - ( 1 - c)^{\delta+1} > 0 \) and \( k > cn \), then the size of the 
largest \(G\)-intersecting, \(k\)-uniform hypergraph is at 
least \( (1 - e^{ -\Omega(n)}) \binom{n}{k} \) (see Theorem~7 of
\cite{BFRT}).  In some sense, this generalizes the trivial observation
that \( \binom{[n]}{k} \) is intersecting for \( k > n/2 \). \\

There is a considerable gap between the values of \(k\) 
for which we have established these two types of  behavior for
maximum \(G\)-intersecting families.  For example, for \(C_n\) we have
(\ref{eq:cycle}) for \( k < C \sqrt n \) while we have 
\( N(C_n,k) > (1 - o(1)) \binom{n}{k} \) for \(k\) greater than
roughly \(.32 n\).  What happens for other values of \(k\)?  Are there
other forms that a maximum \(G\)-intersecting family can take?
Bohman, Frieze, Ruszink\'o and Thoma conjecture that this is not the
case, at least for the cycle.
%
%
\begin{conjecture}
There exists a constant $c$ such that for any fixed $\epsilon>0$
\begin{eqnarray*}
k\leq (c-\epsilon)n & \Rightarrow & N(C_n,k)={n\choose k}-{n-2\choose
k}+{n-4\choose k-2} \\ k\geq (c+\epsilon)n & \Rightarrow &
N(C_n,k)=(1-o(1)){n\choose k} \\ \\
\end{eqnarray*}
\end{conjecture}

The remainder of this note consists of the proof of Theorem~2.

\section{Utilizing $\tau$}
\label{sTAU}

Let $\cal H$ be a hypergraph and $G$ be a graph on vertex set 
$[n]$.  For
$X \subseteq [n]$, we define
\[ N(X):=\{v\in V(G) : v\sim_G w\mbox{ for some }w\in X\}\cup X. \]
For \( x \in [n]\) we write \( N(x) \) for \(N(\{x\})\). 
We will define the hypergraph ${\cal F}$ by setting $f\in {\cal
F}$ if and only if $f=N(h)$ for some $h\in {\cal H}$.  
Note that if $\cal H$ is $G$-intersecting, then
\begin{equation}
\label{eq:cond}
h \in \hee, f \in \fee \Rightarrow h \cap f \neq \emptyset. 
\end{equation}
The quantity
$\tau({\cal F})$ is the cover number of $\cal F$. \\

The proof of Theorem~\ref{tSQRT} follows immediately 
from Lemma~\ref{lemTAU}, which
deals with the case where $\tau_({\cal F})\geq 2$ and
Lemma~\ref{lemONE}, which deals with the case where $\tau({\cal
F})=1$. \\

\begin{lemma} 
Let $G$ be a graph on $n$ vertices with maximum degree $\Delta$ and
clique number \( \omega \), both
constants.  If $ k < \sqrt{ \frac{ \omega n }{2 ( \Delta + 1 )^2 }}$, 
$\cal H$ is a
$k$-uniform, $G$-intersecting hypergraph on $n$ vertices and \(n\) is
sufficiently large, then 
$\tau({\cal F})=1$ or
\begin{equation}
   |{\cal H}|<{n\choose k}-{n-\omega \choose k} \label{eqTAU} .
\end{equation}
\label{lemTAU}
\end{lemma} 

\noindent{\bf Proof.} 

Suppose, by way of contradiction, that
$\tau=\tau({\cal F})\geq 2$ and (\ref{eqTAU}) does not hold.
For \( v \in [n] \) set 
\( \hee_v = \left\{ f \in \fee : u \in f \right\} \), and for 
\( Y \subseteq [n] \) set 
\( \hee_Y = \left\{ f \in \fee : Y \subseteq f \right\} \).  Let
\( \fee_u \) and \( \fee_Y \) be defined analogously \\

We first use \( \tau > 1 \) to get an upper bound \( | \hee_u| \) for 
an
arbitrary \( u \in [n] \).  First 
note that, since \( \tau > 1 \), there exists \( X_1 \in \fee \) such that
\( u \not\in X_1 \).  It follows from (\ref{eq:cond}) that each 
\( f \in \fee_u \) must intersect \(X_1\).  In other words, we have 
\[ \fee_u = \bigcup_{u_1 \in X_1} \fee_{ \{ u,u_1\} }. \]
This observation can be iterated: if
\( i < \tau \) and  
\(Y = \left\{u=u_0,u_1,\ldots,u_{i-1}\right\} \) then there exists 
\( X_{i} \in \fee \) such that \( X_{i} \cap Y = \emptyset \), and 
we have
\[ \fee_{Y} = \bigcup_{u_{i} \in X_{i}} \fee_{ Y \cup \{ u_i\}}. \]
%
%
%
%
Since $|f|\leq (\Delta+1)k$ for all $f \in \fee$, it follows that we have
\begin{equation}
\label{eq:one}
\left|{\cal H}_u\right|\leq 
   \left((\Delta+1)k\right)^{\tau-1}{n-\tau\choose k-\tau}.
\end{equation}\\

On the other hand, by the definition of $\tau$, there exists $v \in [n]$ 
for which
\[ \frac{1}{\tau}\left[{n\choose k}-{n-\omega(G)\choose k}\right]\leq
   \left|{\cal F}_v\right| . \] 
It follows that there exists \(u \in [n]\) such that
\[
\frac{1}{\tau (\Delta+1)}
\left[{n\choose k}-{n-\omega(G)\choose k}\right]\leq
   \left|{\cal H}_u\right| . \]
Applying (\ref{eq:one}) to this vertex we have
\begin{equation*}
   {n\choose k}-{n-\omega(G)\choose k} \leq
   \tau(\Delta+1)^{\tau}k^{\tau-1}{n-\tau\choose k-\tau}.
\end{equation*} \\

In order to show that this is a contradiction, we first note that
$\tau(\Delta+1)^{\tau}k^{\tau-1}{n-\tau\choose k-\tau}$ is a
function that is decreasing in $\tau$.  Indeed, for 
\( \tau \ge 2 \)
we have
%
%
\[ \frac{n-\tau}{k-\tau}\geq\frac{n-2}{k-2}\geq\frac{3}{2}(\Delta+1)k
   \geq\frac{\tau+1}{\tau}(\Delta+1)k  \] 
(note that the condition \( k < \sqrt{ \frac{ \omega n} 
{ 2( \Delta +1)^2}} \) is used in the second inequality).  
It follows
that we have
\[ {n\choose k}-{n-\omega(G)\choose k}\leq 2(\Delta+1)^2k{n-2\choose k-2} , \]
which is not true if $k <\sqrt{\frac{n\omega(G)}{2(\Delta+1)^2}}$
and $n$ is large enough. \hfill $\Box$ \\ \\

\begin{lemma} 
Let $G$ be a graph on $[n]$ with maximum degree $\Delta$, a constant.  
If
$\cal H$ is a $k$-uniform, $G$-intersecting hypergraph on $[n]$, 
$k\leq \sqrt{ \frac{n}{ \Delta (\Delta +1)}}$, $\tau({\cal F})=1$, \(n\)
is sufficiently large
and $\cal H$ is of
maximum size, then there exists a maximum-sized clique $K$ in $G$ such
that $\cal H$ contains every \(k\)-set that intersects $K$.
\label{lemONE}
\end{lemma}

\noindent {\bf Proof.} Let us suppose $\cal H$ is of maximum size and
let $u$ be a cover for $\cal F$, the hypergraph defined above. \\

For $v\in [n]$, let ${\cal H}_v$ denote the members of $\cal
H$ that contain $v$.  Since \( \hee \) is assumed to be 
extremal, we may assume that $|{\cal H}_u|={n-1\choose k-1}$.  
Let \(K\) be the set of $v\in [n]$ such that 
$|{\cal H}_v|={n-1\choose k-1}$.  If $n>(\Delta+2)k$ then 
$K$ must be a
clique in $G$; otherwise, we could find two sets that are not
$G$-intersecting in $\cal H$. \\

We now show that the clique $K$ is maximal.  Assume for the 
sake of contradiction that $v$ is 
adjacent to
every element of \(K\) but 
\(v \not\in K\) (i.e. $|{\cal H}_v|<{n-1\choose k-1}$).
There exists $h \in \cal H$ that \(h\)
contains no member of $N(v)$.  It follows from (\ref{eq:cond}) that we
have
\[
\left|{\cal H}_v\right|<(\Delta+1)k{n-2\choose k-2}. 
\]
Since this bounds holds for all vertices in \( N(u) \setminus K\),
if we have
\begin{equation}
\label{eq:six}
\Delta (\Delta+1)k \binom{n-2}{k-2} < \binom{ n -|K| -1}{k-1}
\end{equation}
then the number of \(k\)-sets that contain \(v\) but do not intersect
\(K\) outnumber those edges in \(\hee \) that contain no 
member of \(K\).  In other words, if (\ref{eq:six}) holds then we get
a contradiction to the maximality of \( \hee \).  However, (\ref{eq:six})
holds for \(n\) sufficiently large (here we use 
\( k < \sqrt{ \frac{ n}{\Delta( \Delta+1)}} \)).

It remains to show that \( K \) is a maximum clique.
Since $K$ is maximal, it must be that any member of $\cal H$ that does
not contain a member of $K$ must contain at least 2 members of
$N(K)\setminus K$.  If
\begin{equation}
   {n\choose k}-{n-|K|\choose k}+{|K|(\Delta- |K|+1)\choose 2}
   {n-|K|-2\choose k-2} < {n\choose k}-{n-|K|-1\choose k}
   \label{eqBOUND}
\end{equation}
and there is some clique of size $|K|+1$, then $\cal H$ cannot be
maximum-sized.  But (\ref{eqBOUND}) holds for $k=o(n)$.  So
the maximum-sized $G$ intersecting family must contain all members of
$\bigcup_{v\in K}{\cal F}_v$ for some $K$ with $|K|=\omega(G)$. \hfill $\Box$
\\ \\

\bibliographystyle{plain}

\begin{thebibliography}{99}
\bibitem{BFRT} T. Bohman, A. Frieze, M. Ruszink\'o, L. Thoma, {\it
$G$-intersecting Families}, Combinatorics, Probability and Computing
{\bf 10}, 376-384.
\bibitem{DF} M. Deza, P. Frankl, {\it Erd\H{o}s-Ko-Rado theorem --
22 years later}, SIAM J. Alg. Disc. Meth. {\bf 4} (1983) 419-431.
\bibitem{EKR} P. Erd\H{o}s, C. Ko, R. Rado, {\it Intersection
theorems for systems of finite sets}, Quart. J. Math. Oxford Ser. 2
{\bf 12} (1961), 313-320.
\end{thebibliography}

\end{document}